\documentclass[10pt]{amsart}

\usepackage{amsfonts}
\usepackage{latexsym}
\usepackage{color}
 \RequirePackage{amsbsy} 
 \RequirePackage{amsopn} 
 \RequirePackage{amsmath}
  \RequirePackage{amssymb}
  \usepackage[mathscr]{eucal}
\definecolor{darkmagenta}{rgb}{0.5, 0, 0.5}
\definecolor{darkblue}{rgb}{0.1, 0.1, 0.7}
\definecolor{darkgreen}{rgb}{0.1, 0.35, 0.1}

 \newtheorem{thee}{Theorem}

 \newtheorem{exxe}[thee]{Example}
 \newtheorem{reem}[thee]{Remark}

 \newcommand{\balf}
 {\renewcommand{\theenumi}{(\alph{enumi})}
 \renewcommand{\labelenumi}{\theenumi}
                      \begin{enumerate}}
\newcommand{\ealf}   {\end{enumerate}
                      \renewcommand{\theenumi}{\arabic{enumi}}
                      \renewcommand{\labelenumi}{\theenumi.}}
\newcommand{\bara}   {\renewcommand{\theenumi}{(\arabic{enumi})}
                      \renewcommand{\labelenumi}{\theenumi}
                      \begin{enumerate} }
\newcommand{\eara}   {\end{enumerate}
                      \renewcommand{\theenumi}{\arabic{enumi}}
                      \renewcommand{\labelenumi}{\theenumi.}}

 \newcommand{\brom}   {\renewcommand{\theenumi}{(\roman{enumi})}
                      \renewcommand{\labelenumi}{\theenumi}
                      \begin{enumerate} }
\newcommand{\erom}   {\end{enumerate}
                      \renewcommand{\theenumi}{\arabic{enumi}}
                      \renewcommand{\labelenumi}{\theenumi.}}

	  \newcommand{\Kr}{\mbox{\rm Kr}}

	   \newcommand{\Max}{\mbox{\rm Max}}

     \newcommand{\f}{\boldsymbol{f}}   
       \newcommand{\F}{\boldsymbol{F}}   
  \newcommand{\co} {\boldsymbol{c}} 
  
  \DeclareMathOperator{\calS}     {\mathcal S}%
  \DeclareMathOperator{\calbV} {\boldsymbol{\mathcal V}}%
  \DeclareMathOperator{\calbW}     {\boldsymbol{\mathcal W}}%
  \DeclareMathOperator{\calJ} {\mathcal J}%
  \RequirePackage{amssymb} 
  \RequirePackage{amsmath} 

\begin{document}

 \title[An $\boldsymbol{e.a.b.}$  not $\boldsymbol{a.b.}$ star operation]{Cancellation properties in ideal systems:\\ an $\boldsymbol{e.a.b.}$  not $\boldsymbol{a.b.}$ star operation}

 \author{Marco Fontana, K. Alan Loper and Ry\^{u}ki Matsuda}
 


\address{M.F.: \ Dipartimento di Matematica, Universit\`a degli Studi ``Roma Tre'', 00146 Rome, Italy.}
\email{fontana@mat.uniroma3.it }
\address{K.A.L.: \ Department of Mathematics, Ohio State University, Newark, Ohio 43055, USA.}
\email{lopera@math.ohio-state.edu}
\address{R.M.: \ Department o.f Mathematics, Ibaraki University, Mito, Ibaraki 310-8512, Japan.} 
\email{rmazda@adagio.ocn.ne.jp}

\date{\today}


 \maketitle

 \begin{abstract} 
 We show that Krull's {\ \texttt{a.b.}}\!   cancellation condition is a properly  stronger condition than Gilmer's {\ \texttt{e.a.b.}}\!   cancellation condition for star operations.
 \end{abstract}
%
 
 
 \section{Introduction}
 

Let $D$ be an integral domain with quotient field $K$. 
Let  $\boldsymbol{{F}}(D) $   [respectively,  $\boldsymbol{{f}}(D)$] be 
the set of all nonzero  fractional ideals [respectively, nonzero finitely 
generated fractional ideals] of $D$.

 A \emph{star operation $\ast$ on $D$} is a mapping $\,\ast : 
\boldsymbol{{F}}(D) \rightarrow
\boldsymbol{{F}}(D)\,$,  $\,E \mapsto E^\ast\,$ such that the following properties hold: $\mathbf{ \bf (\ast_1)}$ $(zD)^\ast = z D$ and $(zE)^\ast = zE^\ast$,\  $\mathbf{ \bf (\ast_2)}$    $E \subseteq F \Rightarrow E^\ast
\subseteq
F^\ast$, \ 
 $\mathbf{ \bf (\ast_3)}$  $E \subseteq 
E^\ast 
 \textrm {   and  
}  E^{\ast \ast} := (E^\ast)^\ast = E^\ast\,,$
 for
all nonzero $\,z \in K\,$, and for all $\,E,F \in
\boldsymbol{{F}}(D)\,$.

\smallskip

Examples of star operations include the \emph{$v$--operation}, defined by $E^{v} : = (D:(D:E))$, for each 
   $E \in
\boldsymbol{{F}}(D)$ \cite[page 396]{G}; the  \emph{$t$--operation}, defined by $E^{t} : = \bigcup\{ F^v\mid F\in \f(D),\  F \subseteq E\}$, for each 
   $E \in\boldsymbol{{F}}(D)$ \cite[page 406]{G};    the {\it $w$--operation} (with the notation proposed by Wang-McCasland)  defined by $E^w := \bigcap \{ED_Q \mid Q \in \Max^t(D) \}$ (where $\Max^t(D)$ is the (nonempy) set of all maximal $t$-ideals of $D$) for all $E \in \F(D)$  \cite{WMc1}.

Let 
     $\ast$ be a star operation on $D$. If $F$ is in $\boldsymbol{f}(D)$, we say that 
     $F$ is 
            \it  $\ast$--\texttt{eab} \rm [respectively,  \it  
             $\ast$--\texttt{ab}\rm],
  if the inclusion
$(FG)^{\ast}
             \subseteq (FH)^{\ast}$ implies that $G^{\ast}\subseteq H^{\ast} $, with $G,\ H \in 
\boldsymbol{f}(D)$, [respectively,  with $G,\ H \in 
 {\F}(D)$].
       
      The operation $\ast$ is said to be   \it \texttt{eab} \rm [respectively, \it \texttt{ab}\rm\ \!]  if each $F\in \boldsymbol{f}(D)$ is $\ast$--\texttt{eab}  [respectively, $\ast$--\texttt{ab}]. \  An \texttt{ab} operation is obviously an \texttt{eab} operation.  Recall also that $E \in \boldsymbol{{F}}(D)$ is called {\it a} ({\it fractional }\!) {\it $\ast$-ideal} of $D$ if $E = E^\ast$.      
      
       In the classical (Krull's) setting, the study of Kronecker function rings on an integral domain generally focusses on the collection of  
{\it ``{\bf\texttt{a}}rithmetisch {\bf\texttt{b}}rauchbar''} (for short, \texttt{a.b.} or, simply, \texttt{ab}, as above) $\ast$--operations \cite{Krull:1936}.   Gilmer's presentation of Kronecker function rings  
 \cite[Section 32]{G} makes use of the (presumably  larger  class of)   {\it ``{\bf\texttt{e}}ndlich {\bf\texttt{a}}rithmetisch {\bf\texttt{b}}rauchbar''}  (for short, \texttt{e.a.b.} or, simply, \texttt{eab}, as above) $\ast$--operations. In this paper, we show that the {\ \texttt{e.a.b.}}\!   cancellation condition is really strictly weaker than the {\ \texttt{a.b.}}\!   cancellation condition. This goal is reached by modifying an example given in the recent paper \cite{FL-2009}.

\section{The Example}

In  \cite[Example 16]{FL-2009}, the authors consider the following example.

  Let $k$ be a field,   $X_1, X_2, ... ,X_n, ...$    an infinite set of indeterminates over $k$ and   $N:=(X_1, X_2, ... ,X_n, ...)k[X_1, X_2, ... ,X_n, ...]$. Clearly, $N$ is a maximal ideal in $k[X_1, X_2, ... ,X_n, ...]$. Set  $D := k[X_1, X_2, ... ,X_n, ...]_N$, let  $M: =ND$  be the maximal ideal of the local domain $D$ and $K:=k(X_1, X_2, ... ,X_n, ...)$   the quotient field of $D$.
 Note that $D$ is a UFD and consider $\boldsymbol{\calbW}$ the set of all the rank one  essential valuation overrings of $D$. Let $\wedge_{\calbW}$ be  the star \texttt{ab} operation on $D$ defined by $\calbW$   \cite[page 398]{G}, {\sl{i.e.}}, for each $E\in \boldsymbol{F}(D)$,
     $$
E^{\wedge_{\calbW}} := \bigcap\{EW \mid W \in \calbW \}.
$$
  It is well known that the $t$--operation  on $D$ is an \texttt{ab} star operation, since $F^t = F^{\wedge_{\calbW}}$  for all  $F \in \f(D)$ \cite[Proposition 44.13]{G}  (more precisely, in this case,
 we have $v=t =w = {\wedge_{\calbW}}$).
  
  Consider the following subset of fractional ideals of $D$:
  $$  {\calJ}:= \{ xF^{t},\ yM,\ zM^2 \mid x, y, z \in K \setminus \{0\}, \;   F\in \boldsymbol{f}(D)\}\,.
  $$
 Since each nonzero principal fractional ideal  of $D$ is in $ {\calJ}$ and, for each ideal $J \in  {\calJ}$ and for each  nonzero $a \in K$, the ideal $aJ$ belongs to ${\calJ}$, then, as above, \cite[Proposition 32.4]{G} guarantees that the   set $ {\calJ}$ defines on $D$ a star operation $\ast$, 
 by setting: 
 $$ E^\ast := \cap \{J \mid J\in \boldsymbol{\calJ}\,,\; J\supseteq E\}\,, \quad  \mbox{for each $E\in \boldsymbol{F}(D)$}\,.
 $$
Since, for each $F \in {\boldsymbol{f}(D)}$,  $F^{t} \in {\calJ}$, it was claimed in \cite[Example 16]{FL-2009} that  $\ast\!\!\mid_{\boldsymbol{f}(D)}= t\!\!\mid_{\boldsymbol{f}(D)}$. This would imply that $\ast$ was an \texttt{eab} operation on $D$, since the operation $t$ \ -- as observed above -- \ is an \texttt{ab} 
star operation on $D$. 

Unfortunately, it is not true that $F^\ast = F^t$ for all $F \in \f(D)$ and, in particular, this equality does not hold if $F \subset D$ and $F^t=D$. For instance, if $I:=(X_1, X_2)$, then clearly, in the Krull domain $D$, we have $I^v= I^t = D$. On the other hand, $I^\ast \subseteq M^\ast = M$, since $M \in {\calJ}$.   More generally, and with a more careful analysis,  we claim that,  if $I := I_{ij}:= (X_i, X_j)$, with $i\neq j \geq 1$, then $I^\ast = M$. 

 {\bf Case 1.} For every $G\in \f(D)$, if $I \subseteq G^t$, then $I \subseteq I^\ast \subseteq M^\ast = M \subsetneq  D = I^t\subseteq G^t$.  Note that the same conclusion holds for every proper ideal $A$ of $D$ such that $A^t =D$, {\sl{i.e.}}, for every $G\in \f(D)$ if $A \subseteq G^t$, then $A \subseteq M^\ast =M \subsetneq G^t$.

 {\bf Case 2.} If   $I \subseteq yM$, for some $0\neq y \in K$, then in particular 
$I \subseteq yD$ and so $D = I^t \subseteq yD$, hence, $y^{-1} \in D$.   There are two possibilities here: either $y^{-1} \in M$ or $y^{-1} \in D\setminus M$. In the first case, {\sl{i.e.}}, if  $y^{-1} \in M$,  then  $1 \in yM$ and so 
$D \subseteq yM$. In the second case, {\sl{i.e.}}, if $y^{-1} \in D\setminus M$, then $y^{-1}$ is invertible in $D$, and so $y, y^{-1} \in D$. Thus, $yM= M$. \\
 Note that the same conclusion holds for every proper ideal $A$ of $D$ such that $A^t =D$, {\sl{i.e.}},
if   $A \subseteq yM$, for some $0\neq y \in K$ and $A^t =D$, then either $D \subseteq yM$ or $ M=yM$.

{\bf Case 3.} If   $I \subseteq zM^2 \subseteq zM$, for some $0\neq z \in K$, then as above $z^{-1} \in D$. Two cases are possible: either $z^{-1} \in M$ or $z^{-1} \in D\setminus M$.
If $z^{-1} \in D\setminus M$, then $z^{-1}$ is invertible in $D$ and so $z, z^{-1} \in D$. Thus, $zM^2= M^2$. However, this is impossible, since $I \not\subseteq M^2$.  If $z^{-1}\in M$, then $M \subseteq zM^2$. \\
 Note that a variation of the previous conclusion holds for every proper ideal $A$ of $D$ such that $A^t =D$ 
  and $A \subseteq M^2$
 (for instance, for $A = I^3$), {\sl{i.e.}},
if   $A \subseteq zM^2$, for some $0\neq z \in K$, $A^t =D$ and $A \subseteq M^2$, then either $A \subseteq zM^2 =M^2$ or $A  \subseteq M^2 \subset M \subseteq zM^2$.

\smallskip

By the previous analysis, we conclude in particular that $I^\ast = \bigcap \{ J\in{\calJ} \mid J \supseteq I \} = M$.
Moreover, since $I^\ast = M$, then we obtain $(I^2)^\ast = (I\cdot I)^\ast = (I^\ast \cdot I^\ast)^\ast = (M^2)^\ast =M^2$.
Furthermore, by the more general analysis for a proper ideal $A$ of $D$ such that $A^t = D$,  in case $A = I^3$  we deduce in particular that $(I^3)^\ast$ also coincides with  $M^2$.
 Therefore,
$$ 
({I}^3)^\ast = M^2 = ({I}^2)^\ast \;\;\;\; \mbox{ but } \;\;\;\; ({I^2})^\ast = M^2 \subsetneq I^\ast  = M\,,$$
and so $\ast$ is not an {\it{\texttt{eab}}} star operation on $D$.

\begin{reem} \rm Let 
  $  {\calJ^{\prime}}:= \{ xD,\ yM,\ zM^2 \mid x, y, z \in K \setminus \{0\}\}\,.
  $
  It is easy to see that \cite[Proposition 32.4]{G} guarantees that the   set $ {\calJ^{\prime}}$ defines on $D$ a star operation that coincides with the star operation $\ast$ defined above by the set $\calJ$, since $F^t = F^v = \bigcap \{ xD\mid x \in K,\, F \subseteq xD\}$, for each $F \in \f(D)$ \cite[Theorem 34.1(1)]{G}. 
  \end{reem}

\smallskip

\bigskip

We provide next  a variation of the previous example in order to construct an \texttt{eab} star operation that   is not  \texttt{ab}.

\begin{exxe} \rm ({\sl Example of an \texttt{eab} star operation that   is not an \texttt{ab} star operation}) 
Let $D$,  $M$  and $K$ be as above.
   Consider the following subset of fractional ideals of $D$:
  $$  {\calS}:= \{ xF^{b},\ yM  \mid x, y  \in K \setminus \{0\}, \;   F\in \boldsymbol{f}(D)\}\,,
  $$
where $b$ is the standard \texttt{ab} operation on $D$ defined by the set  $\calbV$ of all valuation overrings  of $D$, {\sl{i.e.}}, for each $E\in \boldsymbol{{F}}(D)$,
     $$
E^b := E^{\wedge_{\calbV}} := \bigcap\{EV \mid V \in \calbV \}.
$$

 Since each nonzero principal fractional ideal  of $D$ is in $ {\calS}$ and, for each (fractional) ideal $J \in  {\calS}$ and for each  nonzero $a \in K$, the (fractional) ideal $aJ$  belongs to ${\calS}$,  as above, \cite[Proposition 32.4]{G} guarantees that the   set $ {\calS}$ defines on $D$ a star operation $\ast$.

We claim that $\ast$ is an \texttt{eab} operation.  Since the $b$-operation is an \texttt{ab} operation, it is sufficient to prove that 
$\ast\!\!\mid_{\boldsymbol{f}(D)} =  b\!\!\mid_{\boldsymbol{f}(D)}$.  Suppose then that $F \in   {\boldsymbol{f}(D)}$.  Since $F^b \in {\calS}$, it is clear that 
$F^\ast \subseteq F^b$.  Note also that it is well-known that each prime ideal $P$ of an integrally closed domain $D$ is a $b$-ideal, since there always exists a valuation overring of $D$ centered on  $P$ \cite[Theorem 19.6]{G}.  It follows that each ideal of the form $yM$ is a $b$-ideal and, hence, each ideal of 
${\calS}$ is a $b$-ideal.  Since $F^b$ is the intersection of all $b$-ideals which contain $F$, this implies that $F^b \subseteq F^\ast$ (the same conclusion follows also from \cite[Proposition 32.2(b)]{G}).  It follows that 
 $\ast\!\!\mid_{\boldsymbol{f}(D)} =  b\!\!\mid_{\boldsymbol{f}(D)}$ and, hence, $\ast$ is an  \texttt{eab} operation.

  Now, we claim that $\ast$ is not an  \texttt{ab} operation on $D$.
  
To show this, we let
$ I := (X_1, X_2)$  and we prove that $(IM)^\ast = I^\ast = I$.  This will show that $\ast$ is not \texttt{ab}, because we clearly cannot cancel $I$ in the previous equation, {\sl{i.e.}}, $(IM)^\ast =(ID)^\ast$ but   $M^\ast =M \neq D =D^\ast$.

 Therefore, we try to determine which (fractional) ideals in $ {\calS}$ contain $IM$.  We know that $I$ is
 in $ {\calS}$ (since $I \in \f(D)$ and $I$ is a prime ideal of $D$, thus, $I =I^b$) and $I$ contains $IM$.  What we really want to prove is that any (fractional) ideal
 in
 $ {\calS}$ which contains $IM$ also contains $I$.

(1) First, suppose that $IM \subseteq yM$ for some nonzero element $y \in K$.
This causes no problems if  it also implies that $D \subseteq yM$, since then, in particular,
we have $I \subseteq yM$, which is what we want.

 Assume that $y$ is a nonzero element of  $K$ and  that $D \not \subseteq yM$. There are four possibilities here.
 
-- (1, a) If  $y$  is  not in $ D$ and $y^{-1}$ is not in $D$,  then  $yM \cap D
\subseteq yD \cap D  \neq D$. Hence,  $yD \cap D$ is a proper divisorial ideal of $D$ containing $IM$.
This contradicts the fact that $(IM)^v = D$.

-- (1, b) If  $y$  is  not in $D$ and $y^{-1}$ is in $D$, then $y^{-1}$ is in $M$
(since $D$ is local) and so $D \subseteq yM$, which is a contradiction.

-- (1, c) If $y$ is in $D$ and $y$ is invertible in $D$, then $yM =M$, and so
in this case $I \subseteq yM$, which is what we want.

-- (1, d) If $y$ is in $D$ and $y$ is not invertible in $D$, then $IM \subseteq yM \subseteq
yD \subseteq M \neq D$. Again, this contradicts $(IM)^v = D$.

  (2) Now suppose that $G\in \f(D)$ is such that
$IM \subseteq G^\ast = G^b$.
 We extend everything to the $b$-Kronecker function ring of $D$, which is the  following subring of the field of rational functions in one indeterminate, denoted by $T$, over $K$,  {\sl{i.e.}}: 
 $$ \Kr(D, b)\! :=\! \{f/g \in K(T) \mid f,g \in D[T], \ 0\neq g,\ \co(f) \subseteq \co(g)^b \} \!=\! \bigcap \{ V(T) \mid V \in \calbV \},$$
where $\boldsymbol{c}(h)$ is the content of a polynomial $h \in D[X]$ and $V(T):= \{f/g \in K(T) \mid  f, g \in V[T], 0 \neq g \mbox{ and  }  \co(g) =V \}$ is {\it the trivial valuation extension of $V$ to $K(T)$} \cite[definitions at pages 218 and 401, Theorems 32.7 and 32.11, Proposition 33.1]{G}.
 Then,  we should have $I\Kr(D, b)M\Kr(D, b) \subseteq G^b\Kr(D, b) = G\Kr(D, b)$.  Recall that $\Kr(D, b)$ is a B\'ezout domain  and so both $I\Kr(D, b)$ and $G\Kr(D, b)$
 are principal ideals.  This means that we actually have
$M\Kr(D, b)\subseteq G\Kr(D, b)(I\Kr(D, b))^{-1}$, the latter (fractional) ideal being principal.

 There are two possibilities here.

 -- (2, a) $\Kr(D,b) \subseteq G\Kr(D, b)(I\Kr(D, b))^{-1}$.  This would imply that
  $I\Kr(D, b) \subseteq G\Kr(D, b)$.  This would in turn imply that
  $I = I^b  \subseteq G^b = G^\ast$, which was our goal.

-- (2, b) $\Kr(D,b) \not \subseteq G\Kr(D, b)(I\Kr(D, b))^{-1}$.  Rename the principal (fractional) ideal $G\Kr(D, b)(I\Kr(D, b))^{-1}$ as $\mathcal H$.  We know that
$ M\Kr(D, b) \subseteq \mathcal H$. 

 If $\mathcal H$ is an integral ideal of $\Kr(D,b)$, then  obviously
 $ M\Kr(D, b)$ is contained in a proper principal ideal of $\Kr(D,b)$. On the other hand, if $\mathcal H$ is not an integral ideal, then $\mathcal H \cap \Kr(D,b)$ is a proper
integral ideal of $\Kr(D,b)$. Moreover, it is also finitely generated \cite[Proposition 25.4(1)]{G} (hence, principal) in the B\'ezout domain $\Kr(D,b)$. 

Therefore, in either case $ M\Kr(D, b)$  is contained in a proper principal
 ideal
 of $\Kr(D, b)$.  This will lead to a contradiction.
As a matter of fact, suppose  that $\varphi \in \Kr(D,b)$ is a nonzero nonunit rational function and that
$M\Kr(D, b) \subseteq \varphi\Kr(D,b)$.  This means that, for any natural number $n \geq 1$, we have
$X_n \in \varphi\Kr(D, b)$.  On the other hand,  
there are only a finite number of  $X_n$ that are part of
 the reduced representation of $\varphi$.  Without loss of generality, suppose  that these finitely many indices are $1,2,...,r$, {\sl{i.e.}}, $\varphi \in k(X_1, X_2, ..., X_r; T) \ (\subset K(T)) $.  Since $\varphi$ is a nonunit
 in $\Kr(D,b)$,  there must be a valuation overring $V$ of $D$ such that $\varphi$
 is
 a nonunit in the valuation overring  $V(T)$ of $\Kr(D,b)$.  Contract $V$ to the subfield $k(X_1, X_2,...,X_r)$ of $K$.  Call this
valuation domain
 $V_r$.   Then,
 extend
 $V_r$ trivially to $K$.  Call this valuation domain $W$, {\sl{i.e.}}, $W :=V_r(X_{r+1}, X_{r+2}, ....)$.   Clearly, $W$ is a valuation overring of $D$. Then we have a
contradiction,   because $\varphi$ is still a nonunit in the valuation overring $W(T)$ of $\Kr(D, b)$ and each $X_n$ with $n > r$
is
 a unit in $W(T)$.  This contradicts the fact that each $X_n$ lies in
 the
principal ideal $\varphi\Kr(D,b)$.

Therefore, Possibility (2, b) does not occur.  
Therefore, we have to fall back on Possibility (2, a) which
 implies that $I \subseteq 
 G^b = G^\ast$, which was what we needed.
\end{exxe}
\medskip
\noindent
{\bf Acknowledgment.} The first-named author was partially supported by a MIUR-PRIN grant 2008-2011, No. 2008WYH9NY.

\end{document}